\documentclass[11pt]{article}
\usepackage{amsfonts}
\usepackage{amssymb,amsmath,euscript}

\newtheorem{proposition}{Proposition}[section]

\newtheorem{lemma}[proposition]{Lemma}
\newtheorem{theorem}[proposition]{Theorem}

\def\r{\rangle}
\def\dim{{\rm dim}}
\def\dimh{{\rm dim}_{_{\rm H}}}

 \def\beqlb{\begin{eqnarray}}\def\eeqlb{\end{eqnarray}}
 \def\beqnn{\begin{eqnarray*}}\def\eeqnn{\end{eqnarray*}}
 \def\beqn{\begin{equation}}
 \def\eeqn{\end{equation}}
  \def\blemma{\begin{lemma}\sl{}\def\elemma{\end{lemma}}}

\def\R{{\mathbb R}}

\def\d{{\mathbb R}^d}
\def\e{\epsilon}

\def\be{\beta}
\def\la{\lambda}

\def\ga{\gamma}

\def\r{\rho}
\def\ep{\varepsilon}

\def\dim{{\rm dim}_{_{\rm H}}}

\def\E{{\mathbb E}}
\def\P{{\mathbb P}}

\def\N{{\mathbb R}^N}

\def\Var{{\rm Var}}

\def\Gr{{\rm Gr}}

\def\to{\rightarrow}

\def\({\bigg(}
\def\){\bigg)}

\makeatletter \@addtoreset{equation}{section} \makeatother
\newcommand {\qed}%
{%
    {}\hfill
    {}\hfill
    {$\square $}%
    \vspace {0.3cm}%
    \pagebreak [2]%
    \par
}%
\newenvironment{remark}{%
    \vspace{0.3cm} \pagebreak [2]%
    \par%
    \refstepcounter{proposition}
    \noindent%
   {\bf Remark~\theproposition\  }}{\ }%
\begin{document}
%


\title{Spectral Conditions for Strong Local Nondeterminism
and Exact Hausdorff Measure of Ranges of
Gaussian Random Fields }
\author{Nana Luan\\
 University of International Business and Economics, Beijing \\
\and Yimin Xiao \footnote{Research partially
        supported by  NSF grant DMS-1006903.}
\\Michigan State University\\}

\maketitle

\begin{abstract}

Let $X= \{X(t), t \in \R^N\}$ be a
Gaussian random field with values in $\R^d$ defined by
\[
X(t) = \big(X_1(t), \ldots, X_d(t)\big),\qquad t \in \R^N,
\]
where $X_1, \ldots, X_d$ are independent copies of a real-valued,
centered, anisotropic Gaussian random field $X_0$ which has
stationary increments and the property of
strong local nondeterminism. In this paper we determine the exact
Hausdorff measure function for the range $X([0, 1]^N)$.

We also provide a sufficient condition for a
Gaussian random field with stationary increments to be
strongly locally nondeterministic. This condition is
given in terms of the spectral measures of the Gaussian
random fields which may contain either an absolutely
continuous or discrete part. This result strengthens and
extends significantly the related theorems of Berman (1973,
1988), Pitt (1978) and Xiao (2007, 2009), and will
have wider applicability beyond the scope of the present paper.

\end{abstract}

{\sc Running head}: Strong Local Nondeterminism and Hausdorff
Measure of Gaussian Random Fields\\

{\sc 2010 AMS Classification numbers}: 60G15, 60G17,
60G60, 28A80.\\


\section{Introduction}
\label{sec:Intro}

Let $X = \{X(t), t \in \N\}$ be a Gaussian random field with values
in $\R^d$, where
\begin{equation}\label{1.1}
X(t) = \big(X_1(t), \ldots, X_d(t)\big),\qquad t \in \R^N.
\end{equation}
For brevity we call $X$ an $(N, d)$-Gaussian random field.
Sample path properties of $X$ such as the Hausdorff dimensions of
the range $X([0,1]^N)=\{X(t):t\in
[0,1]^N\}$, the graph $\Gr X([0,1]^N)$ $=\{(t,X(t)):t\in [0,1]^N\}$
and the level set $X^{-1}(x)=\{t\in \N: X(t)=x\}$ ($x \in \R^d$)
have been studied by many authors under various assumptions on the
coordinate processes $X_1, \ldots, X_d$. We refer to Adler (1981),
Kahane (1985) and Xiao (2007, 2009) for further information.

In the cases when $X_1, \ldots, X_d$ are independent copies of an
approximately \emph{isotropic} Gaussian random field $X_0$ [a typical
example is fractional Brownian motion], the problems for finding the
exact Hausdorff measure functions for $X([0,1]^N)$, $\Gr X([0,1]^N)$
and $X^{-1}(x)$ have been investigated by Talagrand (1995, 1998),
Xiao (1996, 1997a, 1997b), Baraka and Mountford (2008, 2011).

The main objective of this paper is to study the exact Hausdorff
measure of the range of Gaussian random fields
which are anisotropic in the time-variable. More specifically, we
consider an $(N, d)$-Gaussian random field $X = \{X(t), t \in \N\}$
whose coordinate processes $X_1, \ldots, X_d$ in (\ref{1.1}) are
independent copies of a centered, real-valued Gaussian field $X_0$
with stationary increments and $X_0(0) = 0$ almost surely; and
we assume there exists a constant vector $H=(H_1,\cdots,H_N)\in
(0,1)^N$ such that the following conditions hold:
\begin{itemize}
\item[(C1).]\ There exists a positive constant $c_{_{1,1}} \ge 1$
such that \beqnn c_{_{1,1}}^{-1}\, \r(s,t)^2\le
\E\left(X_0(s)-X_0(t)\right)^2 \le c_{_{1,1}}\, \r(s,t)^2 \quad
\hbox{ for all }\, s, t \in [0, 1]^N, \eeqnn where $\r(s,t)$ is the
metric on $\R^N$defined by \beqnn
\r(s,t)=\sum_{j=1}^N|s_j-t_j|^{H_j},\qquad \forall s,t\in\R^N.
\eeqnn

\item[(C2).]\  There exists a positive constant $c_{_{1,2}}$
such that for all integers $n\geq1$ and all
$u,t^1,\cdots,t^n\in[0,1]^N$, we have \beqnn
\Var(X_0(u)|X_0(t^1),\cdots,X_0(t^n))\geq c_{_{1,2}}\min_{0\leq
k\leq n}\r(u,t^k)^2, \qquad  (t^0=0).
\eeqnn
\end{itemize}

Section 2 below provides a way to construct a large 
class of Gaussian random fields with stationary increments that
satisfy (C1) and (C2). Further examples can be found in Xiao (2009) and 
Luan and Xiao (2010). Under Condition (C1), the $(N, d)$-Gaussian 
random field $X$  has a version which has continuous sample
functions on $[0, 1]^N$ almost surely. Henceforth we will assume
without loss of generality that the Gaussian random field $X$ has
continuous sample paths. When $\{X_0(t), t \in \R^N\}$ satisfies (C2),
we say that $X_0$ has the property of strong local nondeterminism in
metric $\r$ on $[0, 1]^N$.

Xiao (2009) proved that, if Condition (C1) holds, then with probability 1,
\begin{equation}\label{Eq:dim1}
\dimh \,X \big([0,1]^N\big)=\min\bigg\{d;\ \sum_{j=1}^N\frac{1}
{{H}_j}\bigg\},
\end{equation}
where $\sum_{j=1}^0\frac{1}{ {H}_j}:=0$. In the above, $\dimh$
denotes Hausdorff dimension [cf. Kahane (1985) or Falconer (1990)].
Further analytic and fractal properties of Gaussian random fields
which satisfy Conditions (C1) and (C2) have been studied by Xiao
(2009), Bierm\'e \emph{et al.} (2009), Luan and Xiao (2010),
Meerschaert \emph{et al.} (2011) [see also Benassi \emph{et al.}
(1997), Ayache and Xiao (2005), Wu and Xiao (2009, 2011) for related
results].

The first objective of this paper is to refine (\ref{Eq:dim1}) by
determining the exact Hausdorff measure function for
the range $X([0, 1]^N)$.

\begin{theorem}\label{t1.1}
Let $X = \{X(t), t \in \N\}$ be an $(N, d)$-Gaussian random field
with stationary increments defined by (\ref{1.1}), where $X_1, \ldots,
X_d$ are independent copies of a centered, real-valued Gaussian field $X_0$
with stationary increments and $X_0(0) = 0$. We assume that $X_0$ satisfies
Conditions (C1) and (C2). If $d> \sum_{j=1}^N H_j^{-1}$, then we have \beqnn
 0<\varphi_1\textrm{-}m(X([0,1]^N))<
\infty \quad a.s., \eeqnn
where $\varphi_1$ is the function
\beqnn
\varphi_1(r)=r^{\sum_{j=1}^N H_j^{-1}}\log\log\frac{1}{r}
\eeqnn
and  $\varphi_1\textrm{-}m$ is the corresponding Hausdorff
measure.
\end{theorem}

The following remark is concerned with the cases not covered by Theorem
\ref{t1.1}.
\begin{remark}
\begin{itemize}
\item 
If $d< \sum_{j=1}^N H_j^{-1}$, then Theorem 8.2 in Xiao (2009) implies that
$X([0,1]^N)$ a.s. has interior points and hence has positive
$d$-dimensional Lebesgue measure. In this case, Wu and Xiao (2011) showed that 
$X$ has a jointly continuous local time and provides a lower bound for 
the exact Hausdorff measure (in the metric $\rho$) of the level set $X^{-1}(x)$. 
For fractional Brownian motion and some other isotropic Gaussian random fields, 
the exact Hausdorff measure function for $X^{-1}(x)$ has been determined by 
Xiao (1997b) and Baraka and Mountford (2011). However, no such result has 
been established for \emph{anisotropic} Gaussian random fields.  

\item If $d= \sum_{j=1}^N H_j^{-1}$, then $\dim X([0,1]^N)= d$ a.s. The 
problem to determine the exact Hausdorff measure function for $X([0,1]^N)$ in this
``critical case"  is open and is certainly a deeper question.
\end{itemize}
\end{remark}

It will become clear that the proof of Theorem \ref{t1.1} relies
crucially on Condition (C2)--the property of strong local
nondeterminism, which is useful for studying many other sample path
and statistical properties of Gaussian random fields [cf. Xiao
(2009), Xue and Xiao (2011)]. The second objective of this paper is
to provide a rather general condition for a Gaussian random field
with stationary increments to satisfy both Conditions (C1) and (C2).
This condition is given in terms of the spectral measures of the
Gaussian random fields which may contain either an absolutely
continuous or a discrete part. Theorem $\ref{t2.1}$ extends
 the related theorems of Berman (1973, 1988), Pitt
(1978) and Xiao (2007, 2009),  which will have wider applicability
beyond the scope of the present paper. For example, we can apply
this theorem to prove that the solution of a fractional stochastic
heat equation on the circle ${\mathbb S}_1$ [see Tindel, Tudor and
Viens (2004), Nualart and Viens (2009)] has the property of strong
local nondeterminism in the space variable (at fixed time $t$).
Hence fine properties of the sample functions of the solution can
be obtained by using the results in Monrad and Rootz\'{e}n (1995), Xiao (2009),
Luan and Xiao (2010), and Meerschaert, Wang and Xiao (2011).
Similarly, we can show that the spherical fractional Brownian motion
on ${\mathbb S}_1$ introduced by Istas (2005) is also strongly
locally nondeterministic. Both of these processes share local properties
with ordinary fractional Brownian motion with appropriate Hurst
indices. Details of these results will be given elsewhere.

The rest of this paper is organized as follows. Section 2
gives a sufficient condition for a Gaussian random field with
stationary increments to be strongly locally nondeterministic.
Section 3 is concerned with the exact Hausdorff measure function
for the range of $X$. After recalling the definition of Hausdorff
measure and its basic properties, and establishing some estimates,
we prove Theorem \ref{t1.1}.

We end the Introduction with some notation. The inner product of $s,
\, t \in \mathbb{R}^N$ is denoted by  $\langle s, t\rangle$ and the
Euclidean norm of $ t \in \R^N$ is denoted $\| t\|$. Given two
points $s = ( s_1, \ldots, s_N) \in \R^N$ and $t= (t_1, \ldots, t_N)
\in \mathbb{R}^N$, $s\le t$ (resp. $s<t$) means that $s_i\le t_i$
(resp. $s_i < t_i$) for all $1\le i\le N$. When $s\le t$, we use
$[s,t]$ to denote the $N$-dimensional interval (or rectangle) $
[s,t] = \prod_{i=1}^{N}[s_i,t_i].$ For any $T\subseteq \R^N$,
$f(s)\asymp g(s)$ means the ratio $f(s)/g(s)$ is bounded from below
and above by positive and finite constants which are independent of
$s \in T$.

Throughout this paper we will use $c$ to denote an unspecified positive
and finite constant which may not be the same in each occurrence. More specific
constants in Section $i$ are numbered as $ c_{_{i, 1}}, c_{_{i, 2}},
\ldots$.\\

{\bf Acknowledgement} This paper was written while Nana Luan was
visiting Department of Statistics and Probability, Michigan State
University (MSU) with the support of a  grant from China Scholarship
Council (CSC). She thanks MSU for the good working condition and
CSC for the financial support.

The authors thank the referees for their carefully reading of the manuscript 
and their helpful comments.

\section{Spectral condition for strong local nondeterminism of
Gaussian fields with stationary increments}
\label{Sec:Pre}
One of the major difficulties in studying the probabilistic, analytic
or statistical properties of Gaussian random fields is the complexity of
their dependence structures. In many circumstances, the properties
of local nondeterminism can help us to overcome this difficulty so
that many elegant and deep results for Brownian motion can be extended
to Gaussian random fields; see Berman (1973, 1988), Pitt (1978)
and Xiao (2007, 2009) for further information.
Hence, for a given Gaussian random field, it is an interesting question
to determine whether it satisfies certain forms of local nondeterminism.
In this section we provide a general sufficient condition
for a Gaussian random field with stationary increments to satisfy
Conditions (C1) and (C2).

Let $X_0=\{X_0(t),t\in\N\}$ be a real-valued, centered Gaussian
random field with stationary increments and $X_0(0)=0$. We assume
that $X_0$ has  continuous covariance function
$R(s,t)=\E[X(s)X(t)]$. According to Yaglom (1957), $R(s,t)$ can be
represented as
\begin{equation}\label{Eq:covar}
R(s,t)=\int_{\N}\big(e^{i\left< s,\la\right>}-1\big)\big(e^{-i\left<
t,\la\right>}-1\big)\,F(d\la) +\left<s,Mt\right>,
\end{equation}
where $M$ is an $N\times N$ non-negative definite matrix and
$F(d\la)$ is a nonnegative symmetric measure on $\N\backslash\{0\}$
satisfying
\begin{equation}\label{Eq:int}
\int_{\N}\frac{\|\la\|^2}{1+\|\la\|^2}\,F(d\la)<\infty.
\end{equation}
In analogy to the stationary case, the measure $F$ is called the
spectral measure of $X_0$. If $F$ is absolutely continuous with
respect to the Lebesgue measure in $\R^N$, its density $f$ will be
called the spectral density of $X_0$.

It follows from (\ref{Eq:covar}) that $X_0$ has the following
stochastic integral representation:
\begin{equation}\label{Eq:eks}
X_0(t)\stackrel {d}
{=}\int_{\N}\big(e^{i\left<t,\la\right>}-1\big)\,W(d\la) +\left<
Y,t\right>,
\end{equation}
where $\stackrel {d} {=} $ means equality of all finite dimensional
distributions, $Y$ is an $N$-dimensional Gaussian random vector with
mean 0 and covariance matrix $M$, $W(d\la)$ is a centered
complex-valued Gaussian random measure which is independent of $Y$
and satisfies
\[
\E\big(W(A)\overline{W(B)}\big)=F(A\cap B)\quad
 \textup{and}\quad W(-A)=\overline{W(A)}
\]
for all Borel sets $A, B\subseteq \N$ with finite $F$-measure. The 
above properties of $W(d\la)$ ensures that the stochastic integral in
\eqref{Eq:eks} is real-valued. The
spectral measure $F$ is called the control measure of $W$. Since the
linear term $\left< Y,t\right>$ in (\ref{Eq:eks}) will not have any
effect on the problems considered in this paper, we will from now on
assume $Y=0$. This is equivalent to assuming $M=0$ in
(\ref{Eq:covar}). Consequently, for any $h \in \R^N$ we have
\begin{equation}\label{3.4}
\sigma^2(h)\triangleq\E\big(X_0(t+h)-X_0(t)\big)^2=
2\int_{\N}\big(1-\cos\left< h,\la\right>\big)\,F(d\la).
\end{equation}
It is important to note that $\sigma^2(h)$ is a negative definite
function in the sense of I. J. Schoenberg, which is determined by
the spectral measure $F$. See Berg and Forst (1975) for more
information on negative definite functions. If the function
$\sigma^2(h)$ depends only on $\|h\|$, then $X_0$ is called an
isotropic random field. More generally, if
$\sigma^2(h)\asymp\phi(\|h\|)$ in a neighborhood of $h=0$ for some
nonnegative function $\phi$, then $X_0$ is called approximately
isotropic.

Various centered Gaussian random fields with
stationary increments can be constructed by choosing appropriate
spectral measures $F$. For the well known fractional Brownian motion
$B^H= \{B^H(t), t \in \R^N\}$ of Hurst index $H \in (0, 1)$, its
spectral measure has a density function
\begin{equation}\label{Eq:fbmd}
f_H(\la) = c(H, N) \frac 1 {\|\la\|^{2H + {N} }},
\end{equation}
where $c(H, N)>0$ is a normalizing constant such that $\sigma^2(h)=
\|h\|^{2H}$. Since $\sigma^2(h)$ depends on $\|h\|$ only, the
increments of $B^H$ are isotropic and stationary. Examples of
approximately isotropic Gaussian fields with stationary increments
can be found in Xiao (2007).

A typical example of anisotropic Gaussian random field with
stationary increments can be constructed by choosing the spectral
density
\begin{equation}\label{Eq:ftypical}
f(\la) = \frac1 {\Big(\sum_{j=1}^N |\la_j|^{H_j}\Big)^{2 + Q}},
\qquad \forall \la \in \R^N \backslash \{0\},
\end{equation}
where the constants $H_j \in (0, 1)$ for $j = 1, \ldots, N$ and $Q =
\sum_{j=1}^N H_j^{-1}$. This notation will be fixed throughout the
rest of the paper.

It can be verified that $f(\la)$ in (\ref{Eq:ftypical}) satisfies (\ref{Eq:int}) and
the corresponding Gaussian random field $X_0$ has stationary
increments. In the special case when $H_1= \cdots = H_N = H$, 
(\ref{Eq:ftypical}) is very similar to (\ref{Eq:fbmd}). Consequently, 
$X_0$ shares many properties with fractional Brownian motion.

In general, $X_0$ with spectral density (\ref{Eq:ftypical}) is anisotropic 
in the sense that the sample function 
$X_0(t)$ has different geometric and probabilistic characteristics along 
different directions. This gives more flexibility from modeling point of view.
Moreover, $X_0$ is operator-self-similar with exponent $A =(a_{ij})$,
where $a_{ii} = H_i^{-1}$ and $a_{ij} = 0$ if $i \ne j$. The latter
means that for any constant $c > 0$,
\begin{equation}\label{Eq:OPSS2}
\big\{ X_0(c^A\,t),\, t \in \R^N\big\} \stackrel{d}{=}
\big\{c\,X_0(t),\, t \in \R^N\big\},
\end{equation}
where $c^A$ is the linear operator defined by $c^A =
\sum_{n=0}^\infty \frac{(\ln c)^n A^n} {n!}.$ Xiao (2009) proved that
the Gaussian random field $X_0$ satisfies Conditions (C1)
and (C2), and characterized many sample path properties of the
corresponding $(N, d)$-Gaussian field $X$ in
terms of $(H_1, \ldots, H_N)$ explicitly.

We remark that all centered stationary Gaussian random
fields can also be treated using the above framework. In fact, if $Y
= \{Y(t), t \in \R^N\}$ is a centered, real-valued stationary
Gaussian random field, it can be represented as
$Y(t)=\int_{\N} e^{i\left<t,\la\right>}\, W(d\la)$. Thus the random
field $X_0$ defined by
$$
X_0(t)=Y(t)-Y(0)=\int_{\N}\big(e^{i\left<t,\la\right>}-1\big)\,W(d\la),
\quad \forall\ t \in \R^N
$$
is Gaussian with stationary increments and $X_0(0)=0$. Note that the
spectral measure $F$ of $X_0$ in the sense of (\ref{3.4}) is the
same as the spectral measure [in the ordinary sense] of the
stationary random field $Y$.




The main purpose of this section is to prove a sufficient condition
for a general Gaussian random field $X_0$ with stationary increments
to satisfy Conditions (C1) and (C2). In particular, this condition
implies that $X_0$ is strongly locally
nondeterministic in metric $\r$.

To this end we first introduce some notation and state several
lemmas. For any $\la\in\R^N$ and $h>0$, we denote by $C(\la,h)$ the
cube with side-length $2h$ and center $\la$, i.e.,
\beqnn
C(\la,h)=\big\{x\in\R^N:|x_j-\la_j|\leq h,\, j=1,\cdots,N\big\}.
\eeqnn
For any $g\in L^2(\R^N)$, let
$\widehat{g}(\la)=\int_{\R^N}e^{i\left<\la,x\right>}g(x)dx$ be the Fourier
transform of $g$
and let $L^2(C(0,T))$ denote the subspace of $g\in L^2(\R^N)$ whose
support is contained in $C(0,T)$. In the following, Lemma
$\ref{l2.1}$ is Proposition 4 of Pitt (1975). Lemma $\ref{l2.2}$
is taken from Xiao (2007), which is an extension of  a result of
Pitt (1978, p.326).

\begin{lemma}\label{l2.1}
Let $\widetilde{\Delta}(d\la)$ be a positive measure on $\R^N$. If,
for some constant $h > 0$,  $\widetilde{\Delta} (d\la)$ satisfies
\beqn\label{l2.1.1b} 0
<\liminf_{\|\la\|\rightarrow\infty}\widetilde{\Delta}(C(\la,h))\le
\limsup_{\|\la\|\rightarrow\infty}\widetilde{\Delta}(C(\la,h)) <
\infty,
 \eeqn
then, for every $T>0$ satisfying $ThN<\log2$, there exist positive
and finite constants $c_{_{2,2}}$ and $c_{_{2,3}}$ such that
\beqn\label{l2.1.2} c_{_{2,2}}\int_{\R^N}|\widehat{\psi}(\la)|^2d\la
\le \int_{\R^N}|\widehat{\psi}(\la)|^2
 \widetilde{\Delta}(d\la) \le c_{_{2,3}} \int_{\R^N}|\widehat{\psi}(\la)|^2d\la
 \eeqn
for all $\psi \in L^2(C(0,T))$.
\end{lemma}

\begin{lemma}\label{l2.2}
Let $\Delta_1(d\la)$ be a positive measure on $\R^N$ with density
function $\Delta_1(\la)$. If there exist constants $c_{_{2,4}}>0$ and
$\eta>0$ such that
 \beqn\label{l2.2.1}
 \Delta_1(\la)\geq\frac{c_{_{2,4}}}{\|\la\|^{\eta}}\qquad \textrm{for all
 }\la\in\R^N\textrm{ with } \|\la\| \textrm{ large}.
 \eeqn
 Then for any constants $T>0$ and $c_{_{2,5}}$, there exists a
 positive and finite constant $c_{_{2,6}}$ such that for all
 functions $g$ of the form
  \beqn\label{l2.2.2}
  g(\la)=\sum_{j=1}^na_j\left(e^{i\left<s^j,\la\right>}-1\right),
  \eeqn
 where $a_j\in\R$ and $s^j\in C(0,T)$, we have
  \beqnn
  |g(\la)|\leq c_{_{2,6}}\|\la\|\cdot\left(\int_{\R^N}|g(\xi)|^2\
  \Delta_1(\xi)d\xi\right)^{1/2}
  \eeqnn
 for all $\la\in\R^N$ with $\|\la\|\leq c_{_{2,5}}$.
\end{lemma}

Lemma $\ref{l2.3}$ below is an extension of Proposition 8.4 of Pitt
(1978). It allows us to connect the property of strong local
nondeterminism of a Gaussian random field with a general spectral
measure to that of a Gaussian random field with an absolutely
continuous spectral measure, which has been studied in Xiao (2007,
2009).  

\begin{lemma}\label{l2.3}
Let $\Delta_2(d\la)$ be a positive measure on $\R^N$ and suppose that
for some $h>0$, \beqn\label{l2.3.1} 0 <
\liminf_{\|\la\|\rightarrow\infty}\r(0,\la)^{Q+2}\Delta_2(C(\la,h))\le
\limsup_{\|\la\|\rightarrow\infty}\r(0,\la)^{Q+2}\Delta_2(C(\la,h)) <
\infty.
 \eeqn
 Then for any constant $T>0$ with $ThN<\log2$, there exist positive
and finite constants $c_{_{2,7}}$ and $c_{_{2,8}}$ such that
 \beqn
 \label{l2.3.2}
 c_{_{2,7}}\int_{\R^N}\frac{|g(\la)|^2}{(\sum_{j=1}^N|\la_j|^{H_j})^{Q+2}}\,d\la
 \leq\int_{\R^N}|g(\la)|^2\Delta_2(d\la) \le c_{_{2,8}}\int_{\R^N}\frac{|g(\la)|^2}
 {(\sum_{j=1}^N|\la_j|^{H_j})^{Q+2}}\,d\la
 \eeqn
for all $g(\la)$ of the form $(\ref{l2.2.2})$.
\end{lemma}

\noindent{\bf{Proof}.} First we claim that there is a positive
constant $c \le 1 $ such that
\begin{equation}\label{Eq:212}
\begin{split}
c\, \int_{\R^N}\frac{|g(\la)|^2}{\left(\sum_{j=1}^N|\la_j|^{H_j}\right)^{Q+2}}\, d\la
&\leq \int_{\R^N}\frac{|g(\la)|^2}{\left(1+\sum_{j=1}^N|\la_j|^{H_j}\right)^{Q+2}}\, d\la\\
&\qquad \qquad
\le \int_{\R^N}\frac{|g(\la)|^2}{\left(\sum_{j=1}^N|\la_j|^{H_j}\right)^{Q+2}}\, d\la
\end{split}
\end{equation}
for all functions $g$ of the form $(\ref{l2.2.2})$.

Clearly only the first inequality in (\ref{Eq:212}) needs a proof.
For this purpose, we split the first integral in (\ref{Eq:212}) over
$\{\la: \|\la \|\le c_{_{2,5}}\}$ and  $\{\la: \|\la \|
>c_{_{2,5}}\}$ and apply Lemma $\ref{l2.2}$ with
$$\Delta_1(d \la) = \frac{d\la}{\left(1+\sum_{j=1}^N|\la_j|^{H_j}\right)^{Q+2}}
$$
[which satisfies (\ref{l2.2.1})] to derive
\[
\begin{split}
&\int_{\{\|\la\|\le c_{_{2,5}}\}}\frac{|g(\la)|^2}
{\left(\sum_{j=1}^N|\la_j|^{H_j}\right)^{Q+2}}\, d\la\\
&\le c_{_{2,6}}^2\int_{\{\|\la\|\le c_{_{2,5}}\}} \frac{\|\la\|^2}
{\left(\sum_{j=1}^N|\la_j|^{H_j}\right)^{Q+2}}\, d\la\cdot
\int_{\R^N}|g(\xi)|^2\,   \Delta_1(d \xi)\\
&= c_{_{2,9}}\,\int_{\R^N}|g(\xi)|^2\,   \Delta_1(d \xi),
\end{split}
\]
because the first integral in the second line is convergent. It follows from the above
that
\[
\begin{split}
\int_{\R^N}\frac{|g(\la)|^2}{\left(\sum_{j=1}^N|\la_j|^{H_j}\right)^{Q+2}}
\, d\la &\le c_{_{2,9}}\,\int_{\R^N}\frac{|g(\la)|^2}
{\left(1+\sum_{j=1}^N|\la_j|^{H_j}\right)^{Q+2}}\, d\la\\
&\qquad + \int_{\{\la: \|\la \| >c_{_{2,5}}\}}\frac{|g(\la)|^2}
{\left(\sum_{j=1}^N|\la_j|^{H_j}\right)^{Q+2}}\, d\la\\
&\le c_{_{2,10}}\,\int_{\R^N}\frac{|g(\la)|^2}
{\left(1+\sum_{j=1}^N|\la_j|^{H_j}\right)^{Q+2}}\, d\la.
\end{split}
\]
This verifies the first inequality in (\ref{Eq:212}).

Next we take a constant $s>0$ such that $(T+s)hN<\log2$ and denote
$T_1=T+s$. Let $\varphi\in L^2(C(0,s))$ be a function with the
following property 
\beqn\label{l2.3.3} 
c_{_{2,11}}\leq |\widehat{\varphi}(\la)|^2\cdot\left(1+\r(0,\la)\right)^{Q+2}\leq
c_{_{2,12}} 
\eeqn 
for all $\la\in\R^N,$ where $c_{_{2,11}}$ and
$c_{_{2,12}}$ are positive and finite constants. Such a function
$\varphi$ can be constructed as follows. Observe that the function
$\lambda \mapsto \left(1+\r(0,\la)\right)^{-(Q+2)/2}$ is in
$L^2(\R^N)$. Hence it is the Fourier transform of a function $\kappa
\in L^2(\R^N)$. For the constant $s > 0$ chosen above we consider
the function
\[
P_s (t) = \prod_{j=1}^N\left(
        1 - \frac{|t_j|}{s} \right)^+\qquad\text{for all $t\in\R^N$},
\]
where $a^+:=\max(a\,,0)$ for all real numbers $a$. Then the support
of $P_s$ is $C(0,s)$. Recall that the Fourier transform of $P_s$ is
\[
\widehat{P}_s(\xi) := 2^N \prod_{j=1}^N \frac{1-\cos(s \xi_j)}{ s
\xi_j^2} \qquad\text{for all $\xi\in\R^N.$}
\]

Define $\varphi(t) = \kappa(t)P_s (t)$. Then $\varphi\in L^1(C(0,s))\cap
L^2(C(0,s))$ and its Fourier transform is given by
\[
\begin{split}
\widehat{\varphi}(\la) &= \widehat \kappa \star \widehat{P}_s(\la) \\
&=\int_{\R^N} \frac {2^N} {\big(1+\r(0,\la-\xi)\big)^{(Q+2)/2}} \prod_{j=1}^N
\frac{1-\cos(s \xi_j)}{ s \xi_j^2}\, d\xi.
\end{split}
\]
It is clear that $\widehat{\varphi}(\la)>0$ for all $\la \in\R^N$. Writing
\[
\widehat{\varphi}(\la) \cdot\left(1+\r(0,\la)\right)^{(Q+2)/2}
= \int_{\R^N} \frac {2^N \left(1+\r(0,\la)\right)^{(Q+2)/2}} {\big(1+\r(0,\la-\xi)\big)^{(Q+2)/2}} \prod_{j=1}^N
\frac{1-\cos(s \xi_j)}{ s \xi_j^2}\, d\xi
\]
and using the dominated convergence theorem, we see that 
\[
\lim_{\|\lambda\|\to \infty} \widehat{\varphi}(\la) \cdot\left(1+\r(0,\la)\right)^{(Q+2)/2}
= 2^N\int_{\R^N}  \prod_{j=1}^N
\frac{1-\cos(s \xi_j)}{ s \xi_j^2}\, d\xi.
\]
Hence (\ref{l2.3.3}) follows.

Now we continue with the proof of $(\ref{l2.3.2})$. Let
\beqnn
\widehat{\psi}(\la) := g(\la)\widehat{\varphi}(\la)=
\sum_{j=1}^n a_j\big(e^{i\langle s^j,\la\rangle}-1\big)
\widehat{\varphi}(\la),
 \eeqnn
where $s^j \in C(0, T)$ for $j = 1, \ldots, n$. Since 
$\varphi \in L^1(C(0,s))\cap L^2(C(0,s))$, we use the Fourier inversion formula to 
verify that $\psi\in L^2(C(0,T_1))$. Moreover, by (\ref{Eq:212}) and $(\ref{l2.3.3})$,
there is a constant $c \ge 1$ such that
\beqlb\label{l2.3.4}
c^{-1}\int_{\R^N}|g(\la)\widehat{\varphi}(\la)|^2d\la \le
\int_{\R^N}\frac{|g(\la)|^2}{\left(\sum_{j=1}^N|\la_j|^{H_j}\right)^{Q+2}}d\la
 \leq c\int_{\R^N}|g(\la)\widehat{\varphi}(\la)|^2d\la
\eeqlb
for all functions $g$ of the form $(\ref{l2.2.2})$.

Consider the new positive measure $\widetilde{\Delta}(d\la)$ on $\R^N$
defined by
$\widetilde{\Delta}(d\la)=|\widehat{\varphi}(\la)|^{-2}\Delta_2(d\la)$.
It follows from $(\ref{l2.3.1})$ and $(\ref{l2.3.3})$ that
 \beqnn
\liminf_{\|\la\|\rightarrow\infty}\widetilde{\Delta}(C(\la,h))\geq
 c\liminf_{\|\la\|\rightarrow\infty}\r(0,\la)^{Q+2}\Delta_2(C(\la,h))>0
 \eeqnn
 and
\beqnn
\limsup_{\|\la\|\rightarrow\infty}\widetilde{\Delta}(C(\la,h))\leq
 c\limsup_{\|\la\|\rightarrow\infty}\r(0,\la)^{Q+2}\Delta_2(C(\la,h))< \infty
 \eeqnn
Hence the measure $\widetilde{\Delta}(d\la)$ satisfies (\ref{l2.1.1b}).
We apply Lemma $\ref{l2.1}$ to derive that
\[
\begin{split}
c_{_{2,2}}\int_{\R^N}|g(\la)\widehat{\varphi}(\la)|^2d\la
 &\leq\int_{\R^N}|g(\la)\widehat{\varphi}(\la)|^2\widetilde{\Delta}(d\la)\\
 &=\int_{\R^N}|g(\la)|^2\Delta_2(\la)
 \le c_{_{2,3}}\int_{\R^N}|g(\la)\widehat{\varphi}(\la)|^2d\la.
\end{split}
\]
for all functions $g$ of the form $(\ref{l2.2.2})$ provided 
$s^j \in C(0, T)$ for $j = 1, \ldots, n$. This and $(\ref{l2.3.4})$ 
yield $(\ref{l2.3.2})$. \qed

We are ready to prove the main result of this section.

\begin{theorem}\label{t2.1}
Let $\{X_0(t),t\in\R^N\}$ be a real-valued centered Gaussian random
field with stationary increments and $X_0(0) = 0$. If for some constant $h>0$ the
spectral measure $F$ of $X_0$ satisfies \beqn\label{t2.1.2} 0<
\liminf_{\|\la\|\rightarrow\infty}\r(0,\la)^{Q+2}F(C(\la,h))\le
\limsup_{\|\la\|\rightarrow\infty}\r(0,\la)^{Q+2}F(C(\la,h)) <
\infty, \eeqn then for any $T > 0$ such that $ThN < \log 2$, $X_0$
satisfies Conditions (C1) and (C2) on $C(0, T)$.
\end{theorem}

\noindent\textbf{Proof.} First we verify $X_0$ satisfies Condition (C1). For any
$s, t \in C(0, T)$, we apply the stochastic representation of $X_0$ and Lemma
\ref{l2.3} to write
\begin{equation}\label{Eq:wh}
\begin{split}
\E\big(|X_0(s) - X_0(t)|^2\big) &= \int_{\R^N}\big|e^{i\left<s,\la\right>}-
e^{i\left<t,\la\right>}\big|^2\, F(d \la)\\
&\asymp \int_{\R^N} \frac{\big|e^{i\left<s,\la\right>}-
e^{i\left<t,\la\right>}\big|^2}{\big(\sum_{j=1}^N|\la_j|^{H_j}\big)^{Q+2}}\,
d \la.
\end{split}
\end{equation}
Since it has been proved in Xiao (2009) that
\[
\int_{\R^N} \frac{\big|e^{i\left<s,\la\right>}-
e^{i\left<t,\la\right>}\big|^2}{\big(\sum_{j=1}^N|\la_j|^{H_j}\big)^{Q+2}}\, d \la
\asymp \rho(s,t)^2, \qquad \forall s, t \in C(0, T),
\]
we conclude that $X_0$ satisfies  (C1) on $C(0, T)$.

Now we prove that $X_0$ satisfies Condition (C2) on $C(0, T)$.
Denote $r=\min_{0\leq j\leq n}\r(u,t^j)$.
It is sufficient to prove that for all $a_j \in \R$ $(1\leq j\leq n)$ we
have
 \beqn\label{t2.1.1}
 \E\left(\bigg{|}X_0(u)-\sum_{j=1}^na_j X_0(t^j)\bigg{|}^2\right)\geq c_{_{2,10}}r^2
 \eeqn
and $c_{_{2,10}}$ is a positive constant which is independent of $n,$ $a_j$
and the choice of $\{t^j\}$ and $u$. Again by using the stochastic
representation of $X_0$, the left hand side of $(\ref{t2.1.1})$ can
be written as
 \beqnn
 &&\E\left(\bigg{|}X_0(u)-\sum_{j=1}^na_j X_0(t^j)\bigg{|}^2\right)\cr
 &&=\int_{\R^N}\bigg{|}e^{i\left<u,\la\right>}-1
 -\sum_{j=1}^na_j\left(e^{i\left<t^j,\la\right>}-1\right)\bigg{|}^2F(d\la).
 \eeqnn
Note that the function inside the integral is of the form
$({\ref{l2.2.2}})$. We apply Lemma $\ref{l2.3}$ to get
 \beqnn
 &&\int_{\R^N}\bigg{|}e^{i\left<u,\la\right>}-1-
 \sum_{j=1}^na_j\left(e^{i\left<t^j,\la\right>}-1\right)\bigg{|}^2F(d\la)\cr
 &&\geq
 c_{_{2,7}}\int_{\R^N}\bigg{|}e^{i\left<u,\la\right>}-1-
 \sum_{j=1}^na_j\left(e^{i\left<t^j,\la\right>}-1\right)
 \bigg{|}^2\frac{d\la}{\big{(}\sum\limits_{j=1}^N|\la_j|^{H_j}\big{)}^{Q+2}}.
 \eeqnn
However, it has been proved in Theorem 3.2 of Xiao (2009) that the
last integral is bounded from below by $c_{_{2,11}}r^2$, and
$c_{_{2,11}}$ is
a positive constant which is independent of $n,$ $a_j$
and the choice of $\{t^j\}$ and $u$. This proves
$(\ref{t2.1.1})$ and Theorem $\ref{t2.1}.$ \qed

Theorem $\ref{t2.1}$ can be applied directly to Gaussian random
fields with stationary increments and with discrete spectral measure
(or of mixed form $F=F_{ac}+F_{dis}$). It is useful for analyzing
many space-time Gaussian random fields in the literature; see Xue
and Xiao (2011) and the references therein for some examples. In the
following we give an example of Gaussian random field with discrete
spectral measure $F$.

Let $\{\xi_n, n\in \mathbb{Z}^N\}$ and $\{\eta_n, n\in
\mathbb{Z}^N\}$ be two independent sequences of i.~i.~d. $N(0,1)$
random variables, where $\mathbb{Z}$ is the set of  integers.
Let $\{a_n,n\in\mathbb{Z}^N\}$ be a sequence of
real numbers such that
 \beqnn
 \sum_{n\in\mathbb{Z}^N}a_n^2<\infty.
 \eeqnn
Then
 \beqnn
 Y(t)=\sum_{n\in\mathbb{Z}^N}a_n\big(\xi_n\cos\left<n,t\right>
 +\eta_n\sin\left<n,t\right>\big), \qquad
 t\in\R^N
 \eeqnn
is a centered stationary Gaussian random field with covariance
function
\beqnn
\E(Y(t)Y(s))=\sum_{n\in\mathbb{Z}^N}a_n^2\cos\left<n,t-s\right>.
\eeqnn
Hence the spectral measure $F$ of $Y$ is supported on $\mathbb{Z}^N$ with
$F(\{n\})=a_n^2.$ If we choose $\{a_n\}$ such that as $\|n\| \to \infty$,
 \beqnn
 a_n^2 \asymp \frac{1}{\big{(}\sum_{j=1}^Nn_j^{H_j}\big{)}^{Q+2}},
 \eeqnn
then for any fixed constant $h>1$, $F$ satisfies $(\ref{t2.1.2})$. Consider the
Gaussian random field $\{X_0(t), t \in \R^N\}$ defined by $X_0(t) = Y(t) - Y(0)$.
Theorem $\ref{t2.1}$ implies that, for any constant $T > 0$ with $ThN < \log 2$,
$\{X_0(t),t\in\R^N\}$ satisfies Conditions (C1) and (C2) on $C(0,T)$.

Consequently, many sample path properties of $Y$ such as uniform and local moduli
of continuity, Chung's law of the iterated logarithm, existence and joint
continuity of the local times can be derived from the results in Xiao (2009),
Luan and Xiao (2010), and Meerschaert {\it et al.} (2011).


\section{Exact Hausdorff measure function for the range $X([0, 1]^N)$}

In this section, we determine the exact Hausdorff
measure function for the range of an $(N, d)$-Gaussian random field
$X = \{X(t), t \in \N\}$ defined in (\ref{1.1}), where
$X_1, \ldots, X_d$ are independent copies of a real-valued,
centered Gaussian random field $X_0$ with stationary increments,
which satisfies Conditions (C1) and (C2).

First we recall briefly the definition of Hausdorff measure, an
upper density theorem due to Rogers and Taylor (1961) and two useful
inequalities for large and small tails of the supremum of Gaussian processes.
Then we extend a result of Talagrand (1995) to anisotropic Gaussian random fields,
which is applied to derive an upper bound for the $\varphi_1$-Hausdorff measure
of $X([0, 1]^N)$. Finally we prove a law of the iterated logarithm for the
sojourn time of $X$ and derive a lower bound for the $\varphi_1$-Hausdorff measure
of $X([0, 1]^N)$. 

\subsection{Hausdorff measure}
Let $\Phi$ be the class of functions
$\phi:(0,\delta)\rightarrow(0,1)$ which are right continuous,
monotone increasing with $\phi(0_+)=0$ and such that there exists a
finite constant $c_{_{3,1}}> 0$ for which
 \beqnn
 \frac{\phi(2s)}{\phi(s)}\leq c_{_{3,1}},~~~\textrm{for}~0<s<\frac{1}{2}\delta.
 \eeqnn
For $\phi\in\Phi$, the $\phi$-Hausdorff measure of $E\subseteq \d$
is defined by
 \beqnn
\phi\textrm{-}m(E)=\lim_{\epsilon\rightarrow
 0}\inf\left\{\sum_i\phi(2r_i):\, E\subseteq\bigcup_{i=1}^{\infty}
 B(x_i,r_i),r_i<\epsilon\right\},
 \eeqnn
 where $B(x,r)$ denotes the Euclidean open ball of radius $r$ centered at
 $x$. It is known that $\phi\textrm{-}m$
is a metric outer measure and every Borel set in $\d$ is
$\phi\textrm{-}m$ measurable. We say that a function $\phi$ is an
exact Hausdorff measure function for $E$ if
$0<\phi\textrm{-}m(E)<\infty.$ The Hausdorff dimension of $E$ is
defined by
 \beqnn
 \dim E&=&\textrm{inf}\{\alpha>0;s^\alpha\textrm{-}
 m(E)=0\}\cr
 &=&\textrm{sup}\{\alpha>0;s^\alpha\textrm{-}m(E)=\infty\}.
 \eeqnn
We refer to Falconer (1990) for more properties of Hausdorff measure
and Hausdorff dimension.

 The following lemma can be easily derived from the results in
 Rogers and Taylor (1961), which gives a way to get a lower
  bound for $\phi\textrm{-}m(E)$. For any
 Borel measure $\mu$ on $\d$ and $\phi\in\Phi$, the upper
 $\phi$-density of $\mu$ at $x\in\d$ is defined by
  \beqnn
  \overline{D}_\mu^\phi(x)=\limsup_{r\rightarrow
  0}\frac{\mu(B(x,r))}{\phi(2r)}.
  \eeqnn

\blemma \label{l3.1} For a given $\phi\in\Phi$ there exists a
positive constant $c_{_{3,2}}$ such that for any Borel measure $\mu$
on $\d$ and every Borel set $E\subseteq\d$, we have
 \beqnn
\phi\textrm{-}m(E)\geq c_{_{3,2}}\mu(E)\inf_{x\in
E}\{\overline{D}_\mu^\phi(x)\}^{-1}.
 \eeqnn
\elemma

Now we recall some basic facts about Gaussian processes. Consider
a set $S$ and a centered Gaussian process $\{Y(t), t\in S\}$. We provide $S$
with the following canonical pseudo-metric
  \beqnn
  d(s,t)=\|Y(s)-Y(t)\|_2,
  \eeqnn
 where $\|Y\|_2=(\E(Y^2))^{1/2}$. Denote by $N_d(S,\epsilon)$
 the smallest number of open $d$-balls
 of radius $\epsilon$ needed to cover $S$ and let
 $D=\sup\{d(s,t):s,t\in S \}$ be the $d$-diameter of $S$.

The following lemma is well known. It is a consequence of the
Gaussian isoperimetric inequality and Dudley's entropy bound [see
Talagrand (1995)]. \blemma \label{l3.2}There exists a positive
constant $c_{_{3,3}}$ such that for all $u>0$, we have
  \beqnn
 \mathbb{P}\left\{\sup_{s,t\in
 S}|Y(s)-Y(t)|\geq c_{_{3,3}}\bigg{(}u+\int_0^D\sqrt{\log
 N_d(S,\epsilon)}d\epsilon\bigg{)}\right\}\leq\exp\bigg{(}-\frac{u^2}
{D^2}\bigg{)}.
 \eeqnn
\elemma \blemma \label{l3.3} Consider a function $\Psi$ such that
$N_d(S,\epsilon)\leq\Psi(\epsilon)$ for all $\epsilon>0$. Assume
that for some constant $c_{_{3,4}}\ge 1$ and all $\epsilon>0$ we have
 \beqnn
 \Psi(\epsilon)/c_{_{3,4}}\leq\Psi(\frac{\epsilon}{2})\leq c_{_{3,4}}\Psi(\epsilon).
 \eeqnn
 Then
 \beqnn
 \mathbb{P}\left\{\sup_{s,t\in
 S}|Y(s)-Y(t)|\leq u\right\}\geq\exp\big{(}-c_{_{3,5}}\Psi(u)\big{)},
 \eeqnn
 where $c_{_{3,5}}>0$ is a constant depending only on $c_{_{3,4}}$.
\elemma

This was proved in Talagrand (1993). It gives a general
lower bound for the small ball probability of Gaussian processes.

\subsection{Some basic estimates}

Let $X_0= \{X_0(t),\, t \in \R^N\}$ be a centered Gaussian random
field with stationary increments and satisfying Conditions (C1) and
(C2). Without loss of generality, we assume that ${H}_1, \ldots,
{H}_N$ are ordered as
\begin{equation}\label{Eq:Hs2}
0 <  {H}_1 \le  {H}_2 \le \cdots \le  {H}_N < 1.
\end{equation}

In order to solve some dependence problems that are a major
obstacle,  we consider for any given $0<a<b<\infty$ the random field
 \beqnn
 X_0(a,b,t)=\int_{a<\r(0,\la)\le b}\big(e^{i\left<t,\la\right>}-1\big)\,W(d\la),
 \qquad   t \in \R^N.
 \eeqnn
An essential fact is that if $0<a<b<a'<b'<\infty$, then the
Gaussian random fields $\{X_0(a,b,t), t \in \R^N\}$ and
$\{X_0(a',b',t), t \in \R^N\}$ are independent.

Let $X_1(a,b,t),\cdots,X_d(a,b,t)$ be independent copies of
$X_0(a,b,t)$ and let
\beqnn
X(a,b,t)=\left(X_1(a,b,t),\cdots,X_d(a,b,t)\right),\qquad t\in\R^N.
\eeqnn
Then we have the following lemma. For convenience, we write $I = [0, 1]^N$.
\begin{lemma}\label{l3.1'}
Given any $0<a<b$ and $0< \e<r$, we have
 \beqn\label{l3.1'.1}
 \P\left\{\sup_{t \in I: \r(0,t)\leq r}\|X(a,b,t)\|\leq
 \e\right\}\geq\exp\left(-c\left(\frac{r}{\e}\right)^Q\right),
 \eeqn
 where $0<c<\infty $ is an absolute constant.
\end{lemma}
\noindent\textbf{Proof.} It is sufficient to prove $(\ref{l3.1'.1})$
for $X_0(a,b,t)$. Let $S=\{t\in I:\r(0,t)\leq r\}$ and define a
distance $d$ on $S$ by
 \beqnn
 d(s,t)=\big\|X_0(a,b,s)-X_0(a,b,t)\big\|_2.
 \eeqnn
Then (C1) implies $d(s,t)\leq c_{_{1,1}}\sum_{i=1}^N|s_i-t_i|^{H_i}$
for all $s, t \in I$, independent of the choices of $0<a<b$. It follows that
\beqnn
N_d(S,\e)\leq c\,\Big(\frac{r}{\e}\Big)^Q.
\eeqnn
By Lemma $\ref{l3.3}$ we have
\beqnn \P\left\{\sup_{t \in I: \r(0,t)\leq r}|X_0(a,b,t)|\leq
 \e\right\}\geq\exp\left(-c\left(\frac{r}{\e}\right)^Q\right).
\eeqnn This proves Lemma $\ref{l3.1'}$.\qed

The following truncation inequalities are extensions of those in
Lo\'{e}ve (1977, p.209) for $N=1$ and $(3.4)$ and $(3.5)$ in Xiao
(1996) for $N>1$ and $\rho$ being replaced by the Euclidean metric.

\begin{lemma} \label{l4.2}
There exist positive finite constants $c_{_{3,6}}$ and $c_{_{3,7}}$
such that the following hold.
\begin{itemize}
\item[(i)]\, For any $a>0$ and any $t\in\N$ with $\r(0,t)a\leq 1/N$
we have
 \beqn
 \label{l4.2.1}
 \int_{\{\la:\r(0,\la)\leq a\}}\left<t,\la\right>^2F(d\la)\leq c_{_{3,6}}
 \int_{\N}(1-\cos\left<t,\la\right>)F(d\la).
 \eeqn
\item[(ii)]\, For all $a>0$
 \beqn
 \label{l4.2.2}
 \int_{\{\la:\r(0,\la)>a\}}F(d\la)\leq c_{_{3,7}}a^{-2}.
 \eeqn
\end{itemize}
\end{lemma}

\noindent\textbf{Proof.} Notice that when $\r(0,\la)\leq a$, the condition
$\r(0,t)a\leq 1/N$ implies that $|\left<t,\la\right>| < 1$. It
follows that
\beqnn
1-\cos\left<t,\la\right> \geq
\frac{\left<t,\la\right>^2}{2}\(1-\frac{\left<t,\la\right>^2}{12}\)
\geq\frac{11}{24}\left<t,\la\right>^2.
\eeqnn
Then for any $t\in\N$ with $\r(0,t)a\leq 1/N$ we have
 \beqnn
 \int_{\N}(1-\cos\left<t,\la\right>)F(d\la)&\geq&
 \frac{11}{24}\int_{\{\la:|\left<t,\la\right>|\leq1\}}\left<t,\la\right>^2F(d\la)\cr
 &\geq&\frac{11}{24}\int_{\{\la:\r(0,\la)\leq a\}}\left<t,\la\right>^2F(d\la).
 \eeqnn
 That is
 \beqnn
 \int_{\{\la:\r(0,\la)\leq a\}}\left<t,\la\right>^2F(d\la)\leq \frac{24}{11}
 \int_{\N}(1-\cos\left<t,\la\right>)F(d\la).
 \eeqnn

To prove $(\ref{l4.2.2})$, we make the
following two claims:
\begin{itemize}
\item[(a).] For any $u>0$, if $\la_i\neq0$ for $i=1,\dots,N$, then
 \beqnn
 \frac{1}{2^Nu^Q}\int_{\prod_{i=1}^N[-u^{\frac{1}{H_i}},u^{\frac{1}
 {H_i}}]}\cos\left<t,\la\right>dt= \prod_{i=1}^N \frac{\sin\big(
 u^{\frac{1}{H_i}}\la_i\big)}{u^{\frac{1}{H_i}}\la_i}.
 \eeqnn
\item[(b).] For any $u>0$,
 \beqnn
 \int_{\{\la:\r(0,\la)>
 \frac{1}{u}\}}F(d\la)\leq\frac{c}{2^Nu^Q}
 \int_{\prod_{i=1}^N[-u^{\frac{1}{H_i}},u^{\frac{1}{H_i}}]}dt
 \int_{\N}(1-\cos\left<t,\la\right>)F(d\la).
 \eeqnn
\end{itemize}
Claim (a) is obviously true when $N=1$. Suppose it is true for
$N=k$, then for $N=k+1$, we have
 \beqnn
&&\frac{1}{2^{k+1}u^{\frac{1}{H_1}+\cdots+\frac{1}{H_{k+1}}}}
\int_{\prod_{i=1}^k[-u^{\frac{1}{H_i}},u^{\frac{1}{H_i}}]}dt_1\cdots
 dt_k\int_{[-u^{\frac{1}{H_{k+1}}},u^{\frac{1}{H_{k+1}}}]}
 \cos(t_1\la_1+\cdots+t_{k+1}\la_{k+1})dt_{k+1}\cr
 &&=\frac{1}{2^{k}u^{\frac{1}{H_1}+\cdots+\frac{1}{H_{k}}}}
 \int_{\prod_{i=1}^k[-u^{\frac{1}{H_i}},u^{\frac{1}{H_i}}]}dt_1\cdots
 dt_k\cr
 &&\qquad \quad \times \, \frac{\sin(t_1\la_1+\cdots+
 t_k\la_k+u^{\frac{1}{H_{k+1}}}\la_{k+1})-\sin(t_1\la_1+\cdots+
 t_k\la_k-u^{\frac{1}{H_{k+1}}}\la_{k+1})}{2u^{\frac{1}{H_{k+1}}}\la_{k+1}}\cr
 &&=\frac{1}{2^{k}u^{\frac{1}{H_1}+\cdots+\frac{1}{H_{k}}}}
 \int_{\prod_{i=1}^k[-u^{\frac{1}{H_i}},u^{\frac{1}{H_i}}]}
 \cos(t_1\la_1+\cdots + t_k\la_k)dt_1\cdots
 dt_k\frac{\sin
 u^{\frac{1}{H_{k+1}}}\la_{k+1}}{u^{\frac{1}{H_{k+1}}}\la_{k+1}}\cr
 &&=\frac{\sin
 u^{\frac{1}{H_1}}\la_1}{u^{\frac{1}{H_1}}\la_1}\cdots\frac{\sin
 u^{\frac{1}{H_{k+1}}}\la_{k+1}}{u^{\frac{1}{H_{k+1}}}\la_{k+1}}.
 \eeqnn
Hence claim (a) is true for all $N\geq1$.

By Fubini's theorem and claim (a), we have
 \beqnn
 &&\frac{1}{2^Nu^Q}\int_{\prod_{i=1}^N[-u^{\frac{1}{H_i}},
 u^{\frac{1}{H_i}}]}dt\int_{\N}(1-\cos\left<t,\la\right>)F(d\la)\cr
 &&=\int_{\N}\bigg{[}\frac{1}{2^Nu^Q}
 \int_{\prod_{i=1}^N[-u^{\frac{1}{H_i}},u^{\frac{1}{H_i}}]}
 (1-\cos\left<t,\la\right>)dt\bigg{]}F(d\la)\cr
 &&=\int_{\N}\(1-\prod_{i=1}^N
 \frac{\sin  u^{\frac{1}{H_i}}\la_i}{u^{\frac{1}{H_i}}\la_i}\)F(d\la)\cr
 &&\geq\int_{\N\backslash \{\la: |\la_i|
 \leq(Nu)^{-\frac{1}{H_i}},\, \forall i\}}\(1-\prod_{i=1}^N
 \frac{\sin
 u^{\frac{1}{H_i}}\la_i}{u^{\frac{1}{H_i}}\la_i}\)F(d\la)\cr
 &&\geq c\int_{\N\backslash \{\la: |\la_i|
 \leq(Nu)^{-\frac{1}{H_i}},\, \forall i\}}F(d\la)\cr
 &&\ge c\int_{\{\la:\r(0,\la)>\frac{1}{u}\}}F(d\la).
 \eeqnn
Hence claim (b) is verified.

Now we turn to the proof of $(\ref{l4.2.2})$. With claim (b),
$(\ref{3.4})$ and Condition (C1) in hand, we have for $a > 0$,
 \beqnn
 \int_{\{\la:\r(0,\la)>a\}}F(d\la)&\leq& \frac{ca^Q}{2^N}
 \int_{\prod_{i=1}^N [-a^{-\frac{1}{H_i}},a^{-\frac{1}{H_i}}]}dt
 \int_{\N}(1-\cos\left<t,\la\right>)F(d\la)\cr
 &\leq&\frac{ca^Q}{2^N}\int_{\prod_{i=1}^N
 [-a^{-\frac{1}{H_i}},a^{-\frac{1}{H_i}}]}
 \sum_{i=1}^{N}|t_i|^{2H_i}dt\cr
 &\leq&ca^{-2}.
 \eeqnn
 This finishes the proof of Lemma \ref{l4.2}.
\qed

Lemma $\ref{l4.3}$ gives estimates on the small ball probability of the
$(N, d)$-Gaussian random field  $X$ in (\ref{1.1}).
\begin{lemma}
\label{l4.3} There exist constants $c_{_{3,8}}$ and $c_{_{3,9}}$
such that for all $0< \epsilon<r$,
 \beqn
 \label{l4.3.1}
 \exp\(-c_{_{3,8}}\left(\frac{r}{\e}\right)^Q\)\leq
 \P\left\{\sup_{t\in I:\r(0,t)\leq r}
 \|X(t)\|\leq \e\right\}\leq
 \exp\(-c_{_{3,9}}\left(\frac{r}{\e}\right)^Q\).
 \eeqn
\end{lemma}

\noindent\textbf{{Proof.}} Let $S=\{t\in I:\r(0,t)\leq r\}$. It
follows from (C1) that for all $\e\in (0,r)$,
 \beqnn
N_\r(S,\e)\leq
c\prod_{i=1}^N\left(\frac{r}{\e}\right)^{\frac{1}{H_i}}
=c\left(\frac{r}{\e}\right)^Q:=\psi(\e).
 \eeqnn
Clearly $\psi(\e)$ satisfies the condition in Lemma $\ref{l3.3}$.
Hence the lower bound in $(\ref{l4.3.1})$ follows from Lemma
$\ref{l3.3}$.

The proof of the upper bound in $(\ref{l4.3.1})$ is based on
Condition (C2) and a conditioning argument and is similar to the proof of
Theorem 5.1 in Xiao (2009) [see also Monrad and Rootz\'{e}n(1995)]. We
include it for the sake of completeness.
Let $T=\prod_{i=1}^N\big[0, \left(\frac{r}{N}\right)^{\frac{1}{H_i}}\big]$.
Then $T\subseteq S$.  We divide $T$ into
 \beqnn
 \ell:= \prod_{i=1}^N \bigg(\lfloor \big(\frac{r}{N \e}\big)^{\frac{1}{H_i}}
 \rfloor +1\bigg)\ge \left(\frac{r}{N\e}\right)^Q
 \eeqnn
sub-rectangles of side-lengths $\e^{1/H_i}~ (i=1,\ldots,N)$, where
$\lfloor x\rfloor$ is the largest integer no more than $x$. And
denote the lower-left vertices of these rectangles (in any order) by
$t_{k}\ (k=1,\ldots,\ell)$. Then \beqn\label{l4.3.2}
 \P\left\{\sup_{t\in S}
 \|X(t)\|\leq \e\right\}\leq \P\left\{\sup_{1\leq k\leq \ell}
 \|X(t_{k})\|\leq \e\right\}.
 \eeqn
It follows from Condition (C2) that for every $1\leq k\leq \ell$
\beqnn
\Var\left(X_0(t_{k})|X_0(t_{i}):1\leq i\leq k-1\right)\geq c_{_{1,2}}\, \e^{2}.
\eeqnn
By this and Anderson's inequality for Gaussian measures [see Anderson (1995)], we have the
following upper bound for the conditional probabilities
\beqn\label{l4.3.3}
\P\big\{\|X(t_{k})\|\leq\e\big{|}X(t_{i}):1\leq i\leq
k-1\big\}\leq \Phi\Big(\frac 1 {\sqrt{c_{_{1, 2}}}}\Big)^d,
\eeqn
where $\Phi(x)$ is the distribution function of a standard normal
random variable. It follows from $(\ref{l4.3.2})$ and
$(\ref{l4.3.3})$ that
 \beqnn
 \P\left\{\sup_{t\in S}\|X(t)\|\leq
 \e\right\}\leq \Phi\Big(\frac 1 {\sqrt{c_{_{1, 2}}}}\Big)^{\ell d}\le
 \exp\(-c_{_{3,9}}\left(\frac{r}{\e}\right)^Q\).
 \eeqnn
 Thus we obtain the upper bound in $(\ref{l4.3.1})$.\qed

 The main estimate is given in the following proposition.

 \begin{proposition}
 \label{p4.4}
 There exist positive constants $\delta_1$ and $c_{_{3,10}}$ such that for any
 $0<r_0\leq \delta_1$, we have
 \beqn\label{Eq:4.4}
 \begin{split}
 &\P\left\{\exists \,r\in[r_0^2,r_0], \sup\limits_{t \in I: \r(0,t)\leq r}\|X(t)\|\leq
 c_{_{3,10}} r\left(\log\log\frac{1}{r}\right)^{-1/Q}\right\}\\
 &\qquad \qquad \geq 1-\exp\left(-\left(\log\frac{1}{r_0}\right)^{1/2}\right).
 \end{split}
 \eeqn
 \end{proposition}

\noindent\textbf{Proof.} Though the main idea of the proof is similar to
the proof of Proposition 4.1 in Talagrand (1995), some modifications
are needed to characterize the anisotropic nature of $X$. Let $U>1$ be
a number whose value will be determined later. For $k\geq0$, let
$r_k=r_0U^{-2k}$.  Consider the largest integer $k_0$ such that
 \beqlb
 k_0\leq\frac{\log(1/r_0)}{2\log U}.
 \eeqlb
Thus, for $k\leq k_0$ we have $ r_0^2\leq r_k\leq r_0$. It thereby
suffices
 to prove that
\beqnn
 \P\left\{\exists k\leq k_0, \sup\limits_{t \in I: \r(0,t)\leq r_k}\|X(t)\|\leq
 c\, r_k\left(\log\log\frac{1}{r_k}\right)^{-1/Q}\right\}
 \geq 1-\exp\left(-\left(\log\frac{1}{r_0}\right)^{1/2}\right).
 \eeqnn
Let $a_k=r_0^{-1}U^{2k-1}$ and we define for $k=0,1,\cdots$
 \beqnn
 X_{0,k}(t)=X_0(a_k,a_{k+1},t)
 \eeqnn
 and
 \beqnn
\widehat{X}_k(t)=\big(X_{1,k}(t),\cdots, X_{d,k}(t)\big),
 \eeqnn
where $X_{1,k}(t),\cdots, X_{d,k}(t)$ are independent copies of $
 X_{0,k}(t)$. Furthermore, we assume $X_1-X_{1,k},\cdots,X_{d}-X_{d,k}$ are independent copies of
$X_0-X_{0,k}$. We note that the Gaussian random fields
$\widehat{X}_0,\widehat{X}_1,\cdots$ are independent. By Lemma
$\ref{l3.1'}$ we can find a constant $c_{_{3,11}} >0$ such that, if
$r_0$ is small enough, then for each $k\geq0$
 \beqlb
 &&\P\left\{\sup\limits_{t \in I: \r(0,t)\leq r_k} \|\widehat{X}_k(t)\|
 \leq c_{_{3,11}}\, r_k\left(\log\log\frac{1}{r_k}\right)^{-1/Q}\right\}\cr
 &&\geq \exp\left(-\frac{1}{4}\log\log\frac{1}{r_k}\right)=\frac{1}
 {(\log1/r_k)^{\frac{1}{4}}}\cr
 &&\geq\frac{1}{(2\log1/r_0)^{\frac{1}{4}}}.
 \eeqlb
By independence,
 \beqlb
 \label{p4.4.1}
 &&\P\left\{\exists k\leq k_0,\sup\limits_{t \in I: \r(0,t)\leq r_k}
 \|\widehat{X}_k(t)\|\leq
 c_{_{3,11}}r_k\left(\log\log\frac{1}{r_k}\right)^{-1/Q}\right\}\cr
 &&\qquad\qquad\geq1-\left(1-\frac{1}{(2\log1/r_0)^{1/4}}\right)^{k_0}\cr
 &&\qquad\qquad\geq
 1-\exp\left(-\frac{k_0}{(2\log1/r_0)^{1/4}}\right),
 \eeqlb
where the last inequality follows from the elementary inequality $1 - x \le e^{-x}$
for all $x \ge 0$.

Let $\be=\min\{\frac{1}{H_N}-1,2\}$. We claim that for any $u\geq
cr_kU^{-\frac{\be}{2}}\sqrt{\log U}$,
 \beqn
 \label{p4.4.2}
 \P\left\{\sup\limits_{t \in I: \r(0,t)\leq r_k} \|X(t)-\widehat{X}_k(t)\|\geq
 u\right\}\leq\exp\left(-\frac{u^2}{cr_k^2U^{-\be}}\right).
 \eeqn
To see this, it's enough to prove that $(\ref{p4.4.2})$ holds for
$X_0- X_{0,k}.$ Consider $S=\{t\in I:\r(0,t)\leq r_k\}$ and on $S$ the
distance
 \beqnn
 d(s,t)=\big\|(X_0(s)-X_{0,k}(s))-(X_0(t)-X_{0,k}(t))\big\|_2.
 \eeqnn
Then $d(s,t)\leq c\sum_{i=1}^N|s_i-t_i|^{H_i}$ and $N_d(S,\e)\leq
c(\frac{r_k}{\e})^Q$. Now we estimate the diameter $D$ of $S$. For
any $t\in S$,
 \beqlb\label{p4.4.4}
 &&\E\left[\left(X_0(t)-X_{0,k}(t)\right)^2\right]
 =2\int_{\{\la:\r(0,\la)\leq a_k\}\cup\{\la:\r(0,\la)>a_{k+1}\}}
 (1-\cos\left<t,\la\right>)F(d\la)\cr
 &&\leq2\int_{\{\la:\r(0,\la)\leq a_k\}}(1-\cos\left<t,\la\right>)F(d\la)
 +4\int_{\{\la:\r(0,\la)>a_{k+1}\}}F(d\la)\cr
 &&=:I_1+I_2.
 \eeqlb
The second term is easy to estimate: By Lemma $\ref{l4.2}$,
 \beqn\label{p4.4.5}
I_2\leq ca_{k+1}^{-2}.
 \eeqn
 For the first term $I_1$, we use the elementary inequality
 $1-\cos\left<t,\la\right>\leq\frac{1}{2}
 \left<t,\la\right>^2$ to derive that for all $t\in S$
 \beqnn
 I_1&\leq&\int_{\{\la:\r(0,\la)\leq a_k\}}\left<t,\la\right>^2F(d\la)\cr
 &=&N^{\frac{2}{H_1}}U^{-\frac{1}{H_N}}\int_{\{\la:\r(0,\la)\leq a_k\}}
 \big<\frac{U^{\frac{1}{2H_N}}}{N^{\frac{1}{H_1}}}t,\la\big>^2F(d\la)\cr
 &=&cU^{-\frac{1}{H_N}}\int_{\{\la:\r(0,\la)\leq a_k\}}\left<t',\la\right>^2F(d\la),
 \eeqnn
where $t'= U^{\frac{1}{2H_N}}N^{-\frac{1}{H_1}}\,t.$
Since
$\r(0,t')\leq\frac{1}{N}U^{\frac{1}{2}}\r(0,t)\leq\frac{1}{N}U^{\frac{1}{2}}
r_k < \frac{1}{Na_k},$
it follows from Lemma $\ref{l4.2}$ and (C1) that
 \beqn
 \label{p4.4.6}
 I_1\leq cU^{-\frac{1}{H_N}}\r(0,t')^2\leq
 cU^{1-\frac{1}{H_N}}\r(0,t)^2\leq cr_k^2U^{-(\frac{1}{H_N}-1)}.
 \eeqn
With $(\ref{p4.4.4})$, $(\ref{p4.4.5})$ and $(\ref{p4.4.6})$ in
hand, the diameter of $S$ satisfies
 \beqlb
 \label{p4.4.7}
 D^2&\leq&c\left[r_k^2U^{-(\frac{1}{H_N}-1)}+a_{k+1}^{-2}\right]\cr
 &\leq&cr_k^2\left[U^{-(\frac{1}{H_N}-1)}+U^{-2}\right]\cr
 &\leq&cr_k^2U^{-\be},
 \eeqlb
 where $\be=\min\{\frac{1}{H_N}-1,2\}$.
Some simple calculations yield
 \beqlb
 \label{p4.4.8}
 \int_0^D\sqrt{\log N_d(S,\e)}d\e
 &\leq& c\int_0^{cr_kU^{-\frac{\be}{2}}}\sqrt{\log
 \frac{r_k}{\e}}\,d\e\cr
 &\leq&cr_kU^{-\frac{\be}{2}}\sqrt{\log U}.
 \eeqlb
Hence we use Lemma $\ref{l3.2}$ and $(\ref{p4.4.8})$ to derive that
for any $u\geq cr_kU^{-\frac{\be}{2}}\sqrt{\log U}$,\beqn
 \label{p4.4.3}
 \P\left\{\sup\limits_{\r(0,t)\leq r_k}|X_0(t)-X_{0,k}(t)|\geq
 u\right\}\leq\exp\left(-\frac{u^2}{cr_k^2U^{-\be}}\right).
 \eeqn
 Thus we have proved (\ref{p4.4.2}).

Now we continue our proof of (\ref{Eq:4.4}). Let $U=(\log1/r_0)^{1/\be}$.
We see that for $r_0 >0$ small
 \beqnn
 U^{\be/2}\left(\log U\right)^{-1/2}\geq
 \left(\log\log\frac{1}{r_0}\right)^{1/Q}.
 \eeqnn
Take $u=c_{_{3,11}}r_k(\log\log1/r_0)^{-1/Q}$. It follows from
$(\ref{p4.4.2})$ that
\beqnn
&&\P\left\{\sup\limits_{t \in I: \r(0,t)\leq r_k}\|X(t)-\widehat{X}_k(t)\|
\geq c_{_{3,11}}r_k\left(\log\log\frac{1}{r_0}\right)^{-1/Q}
\right\}\cr
&&\qquad\qquad\leq\exp\left(-\frac{U^\be}{c_{_{3,12}}
\left(\log\log1/r_0\right)^{2/Q}}\right).
\eeqnn
Combining this with $(\ref{p4.4.1})$, we get
\beqlb\label{p4.4.9}
&&\P\left\{\exists k\leq k_0,\sup\limits_{\r(0,t)\leq r_k}\|X(t)\|\leq
2c_{_{3,11}}r_k\left(\log\log\frac{1}{r_k}\right)^{-1/Q}\right\}\cr
&&\qquad\quad\geq1-\exp\left(-\frac{k_0}{(2\log1/r_0)^{1/4}}\right)\cr
&&\qquad\qquad\quad-k_0\exp\left(-\frac{U^\be}{c_{_{3,12}}
\left(\log\log1/r_0\right)^{2/Q}}\right).
\eeqlb

We recall that
\beqnn
\frac{\log\big(1/r_0\big)}{4\log U}\leq k_0\leq\log\frac{1}{r_0}.
\eeqnn
Then the right-hand side of $(\ref{p4.4.9})$ is at least
$1-\exp(-(\log1/r_0)^{1/2})$ when $r_0 >0$ is small enough. This
completes the proof.
 \qed

\subsection{Upper bound for the Hausdorff measure of the range}

We start with the following result on the uniform modulus of continuity
of $X_0$. See, e.g., Xiao (2009). More precise result can be found in
Meerschaert \emph{et al.} (2011).

\begin{lemma}\label{l4.1}
Let $X_0=\{X_0(t),t\in \N\}$ be a centered Gaussian random field
with values in $\R$. If Condition  (C1) is satisfied, then there
exists a positive and finite constant $c_{_{3,13}}$ such that
\begin{equation}
\limsup_{\|\ep\|\to 0}\frac{\sup_{t\in [0, 1]^N,\,s\in [0,\ep]}
|X_0(t+s)-X_0(t)|} {\rho(0,
\ep)\,\sqrt{\log(1+\rho(0,\ep)^{-1})}}\leq c_{_{3,13}}, \quad \hbox{a.s.}
\end{equation}
\end{lemma}

Now we derive an upper bound for the Hausdorff measure of $X([0,1]^N)$.
\begin{theorem}\label{t4.1}
If $d>Q$, then there exists a constant $c_{_{3,14}} >0$  such that
\begin{equation}
\varphi_1\textrm{-}m(X([0,1]^N)) \le c_{_{3,14}} \qquad a.s.,
\end{equation}
where $\varphi_1(r)=r^Q\log\log1/r$.
\end{theorem}

 \noindent\textbf{Proof.} For $k\geq1$, consider the set
 \beqlb
 R_k&=&\bigg{\{}t\in[0,1]^N:\exists~r\in[2^{-2k},2^{-k}]
 \textrm{ such that}\cr
 &&~~~~~~~~~~~\sup_{s \in I: \r(s,t)\leq r}\|X(s)-X(t)\|\leq
 c_{_{3,10}}\, r(\log\log\frac{1}{r})^{-1/Q}\bigg{\}}.
 \eeqlb

By Proposition $\ref{p4.4}$ we have
 \beqnn
 \P\{t\in R_k\}\geq1-\exp(-\sqrt {k/2}).
 \eeqnn
Denote by $L_N$ the Lebesgue measure on $\R^N$. It follows from
Fubini's theorem that $\P(\Omega_0)=1$, where
 \beqnn
\Omega_0=\left\{\omega:\ L_N(R_k)\geq 1-\exp(-\sqrt{k}/4)\textrm{
infinitely often}\right\}.
 \eeqnn
On the other hand, by Lemma $\ref{l4.1}$, there exists an event
$\Omega_1$ such that $\P(\Omega_1)=1$ and for all
$\omega\in\Omega_1$, there exists $n_1=n_1(\omega)$ large enough
such that for all $n\geq n_1$ and any rectangle $I_n$ of
side-lengths $2^{-n/H_i}(i=1,\cdots,N)$ that meets $[0,1]^N$, we
have
\beqlb \label{Eq:323}
\sup_{s,t\in  I_n}\|X(t)-X(s)\| \leq
c2^{-n}\sqrt{\log[1+(N2^{-n})^{-1}]}\leq
 c2^{-n}\sqrt n.
\eeqlb

Now for a fixed $\omega\in \Omega_0\cap\Omega_1$, we show that
$\varphi_1\textrm{-}m(X([0,1]^N))\le c_{_{3,14}}<\infty$. Consider $k\geq1$ such
that
 \beqnn
 L_N(R_k)\geq 1-\exp(-\sqrt{k}/4).
 \eeqnn
For any $n\geq1,$ we divide $[0,1]^N$ into $2^{nQ}$ disjoint (half-open and half closed)
rectangles of side-lengths $2^{-n/H_i}(i=1,\cdots,N)$. Denote by $I_n(x)$ the
rectangle of side-lengths $2^{-n/H_i}(i=1,\cdots,N)$ containing $x$.
For any $x\in R_k$ we can find the smallest integer $n$ with $k\leq n \leq2k+\ell_0$ (where
$\ell_0$ depends on $N$ only) such that
 \beqlb
 \label{3.52}
 \sup_{s,t\in I_n(x)}\|X(t)-X(s)\|\leq
 c2^{-n}(\log\log2^n)^{-1/Q}.
 \eeqlb
Thus we have
 \beqnn
R_k\subseteq V=\bigcup_{n=k}^{2k+\ell_0}V_n
 \eeqnn
and each $V_n$ is a union of  rectangles $I_n(x)$ satisfying
$(\ref{3.52})$. Clearly $X(I_n(x))$ can be covered by a ball of
radius
 \beqnn
 \rho_n=  c2^{-n}(\log\log2^n)^{-1/Q}.
 \eeqnn
Since $\varphi_1(2\rho_n)\leq c2^{-nQ}=cL_N(I_n)$, we obtain
 \beqlb
 \label{3.19*}
 \sum_{n= k}^{k+\ell_0} \sum_{I_n\in V_n}\varphi_1(2\rho_n)\leq \sum_{n}\sum_{I_n\in
 V_n}cL_N(I_n)
 =cL_N(V)\leq c.
 \eeqlb
Thus $X(V)$ is contained in the union of a
family of balls $B_n$ of radius $\rho_n$ with $\sum_n\varphi_1(2\rho_n)\leq c.$

On the other hand, $[0,1]^N\backslash V$ is contained in a union of
rectangles of side-lengths $2^{-q/H_i}(i=1,\cdots,N)$ where
$q=2k+\ell_0$, none of which meets $R_k$. There can be at most
 \beqnn
 2^{Qq}L_N([0,1]^N\backslash V)\leq c2^{Qq}\exp(-\sqrt k/4)
 \eeqnn
such rectangles. Since $\omega\in\Omega_1$, \eqref{Eq:323} implies that, for each of these
rectangles $I_q$, $X(I_q)$ is contained in a ball of radius
$c2^{-q}\sqrt q$. Thus $X([0,1]^N\backslash V)$ can be covered by a
family $B_n$ of balls of
 radius $\rho_n= c2^{-q}\sqrt q$ such that
 \beqlb
 \label{3.20*}
 \sum_n\varphi_1(2\rho_n)\leq (c2^{Qq}\exp(-\sqrt
 k/4))\cdot(c2^{-qQ}q^{Q/2}\log\log(c2^qq^{-1/2}))\leq1 \eeqlb
 for $k$ large enough. Since $k$ can be arbitrarily large, Theorem
 $\ref{t4.1}$ follows from $(\ref{3.19*})$ and $(\ref{3.20*})$.
 \qed

\subsection{Lower bound for the Hausdorff measure of the range}
\label{sec:lower bounds}

\begin{theorem}\label{t5.1}
If $d>Q$, then there exists a constant $c_{_{3,15}} >0$  such that
\begin{equation}
\varphi_1\textrm{-}m(X([0,1]^N)) \ge c_{_{3,15}} \qquad a.s.,
\end{equation}
where  $\varphi_1(r) = r^Q \log \log 1/r$.
\end{theorem}

In order to prove Theorem \ref{t5.1}, we first study the asymptotic
behavior of the sojourn measure of $X$. For any $r >0$ and $y\in\R^d$,
define
\[
T_{y}(r) = \int_{I} {\bf 1}_{\{\|X(t) - y\|\le r\}}\, dt,
\]
the sojourn time of $X$ in the ball $B(y, r)$. If $y=0,$ we write
$T(r)$ for $T_0(r)$.

\begin{lemma}\label{l5.1}
If $d>Q,$ then there is a finite constant $c_{_{3,16}}$ such that
\beqn\label{l5.1.1}
\E\left( T(r)^n\right)\le c_{_{3,16}}^n n!\,r^{Qn} \eeqn
for all for all integers $n \ge 1$ and $0<r<1$.
\end{lemma}

\noindent\textbf{Proof.} For $n=1$, by Fubini's theorem and (C1) we
have
 \beqnn
 \E \big(T(r)\big)&=&\int_I\P\left\{\|X(t)\|<r\right\}dt\cr
 &\leq&\int_I\min\bigg\{1,c\bigg(\frac{r}{\r(0,t)}\bigg)^d\bigg\}dt\cr
 &=&\int_{\{t\in I:\r(0,t)\leq
 cr\}}dt+c\int_{\{t\in I:\r(0,t)>cr\}}\bigg(\frac{r}{\r(0,t)}\bigg)^d
 dt\cr
 &=:& J_1+J_2.
 \eeqnn
The first term is easy to estimate:
\beqn
\label{t5.1.1}
J_1\leq c\prod_{i=1}^Nr^\frac{1}{H_i}=cr^Q.
\eeqn
For the second term, we use the following elementary fact: Given
positive constants $\be$ and $\gamma$,
there exists a finite constant $c_{_{3,17}}$ such that for all
$a>0,$
\beqlb\label{t5.1.2}
 \int_0^\infty\frac{dx}{(a+x^\be)^\ga}=\begin{cases}
 c_{_{3,17}}a^{-(\ga-\frac{1}{\be})} &\textrm{if~~} \be\ga>1,\\
 +\infty &\textrm{if~~} \be\ga\leq1.
 \end{cases}
 \eeqlb
Since $\r(0,t)>cr$ implies that $t_{j_0}\geq cr^{1/H_{j_{_{0}}}}$
for some $j_{_{0}}\in\{1,\cdots,N\}$, without loss of generality we
assume $j_{_{0}}=1$. Then using $(\ref{t5.1.2})$ $(N-1)$ times, we
obtain \beqlb \label{t5.1.3} J_2&\leq&
cr^d\int_{cr^{\frac{1}{H_1}}}^1dt_1\int_{[0,1]^{N-1}}
\frac{dt_2,\cdots,dt_N}{(\sum_{i=1}^Nt_i^{H_i})^d}\cr
&\leq&cr^d\int_{cr^{\frac{1}{H_1}}}^1dt_1\int_{[0,1]^{N-2}}
\frac{dt_2,\cdots,dt_{N-1}}
{(\sum_{i=1}^{N-1}t_i^{H_i})^{d-\frac{1}{H_N}}}\cr
&\leq&cr^d\int_{cr^{\frac{1}{H_1}}}^1\frac{dt_1}
{(t_1^{H_1})^{d-\sum_{i=2}^N\frac{1}{H_i}}}\cr &\leq& c\, r^Q,
\eeqlb where the last step follows from the assumption that $d>Q$.
It follows from $(\ref{t5.1.1})$ and $(\ref{t5.1.3})$ that
 \beqn\label{t5.1.4}
 \E \big(T(r)\big)\leq cr^Q.
 \eeqn
For $n\geq2,$
 \beqn
 \label{t5.1.5}
 \E (T(r)^n)=\int_{I^n}\P \{\|X(t^j)\|<r,1\leq j\leq n\}dt^1\cdots
 dt^n.
 \eeqn
Consider $t^1,\cdots,t^n\in I$ satisfying
 \beqnn
t^j\neq0,~~\textrm{for }j=1,\cdots,n~~\textrm{and}~~t^j\neq
t^k~~\textrm{for}~~j\neq k.
 \eeqnn
By Condition (C2), we have
\beqlb\label{t5.1.6}
\Var\big(X_0(t^n)\big{|}X_0(t^1),\cdots,X_0(t^{n-1})\big)\geq
c_{_{1,2}} \min\limits_{0\leq k\leq
n-1}\r(t^n,t^k)^2,
\eeqlb
where $t^0=0$. Since conditional distributions in Gaussian
processes are still Gaussian, $(\ref{t5.1.6})$ and Anderson's
inequality yield that for all $x^1, \ldots, x^{n-1}\in \R^d$,
\beqlb
\label{t5.1.7}
\begin{split}
&\P\bigg{\{}\|X(t^n)\|<r\big{|}X(t^1)=x^1,\cdots,X(t^{n-1})
=x^{n-1}\bigg{\}}\\
&~~~~~~\leq  c\min\bigg\{1,\, \bigg(\frac{r} {\min\limits_{0\leq k\leq
n-1}\r(t^n,t^k)}\bigg)^{d}\bigg\}.
\end{split}
 \eeqlb
It follows from (\ref{t5.1.7}) and an argument similar to the proof of
(\ref{t5.1.4}) that
\beqlb \label{t5.1.7b}
\begin{split}
&\int_I \P\bigg{\{}\|X(t^n)\|<r\big{|}X(t^1)=x^1,\cdots,X(t^{n-1})
=x^{n-1}\bigg{\}}dt^n \cr
&~~~~~~ \leq c\int_I\sum_{k=0}^{n-1}\min\bigg\{1,c\bigg(\frac{r}
 {\r(t^n,t^k)}\bigg)^d\bigg\}\,dt^n\cr
&~~~~~~\leq c\,n\int_I\min\bigg\{1,c\bigg(\frac{r}{\r(0,t^n)}\bigg)^d\bigg\}
 \,dt^n\cr
&~~~~~~\leq c\,nr^Q.
 \end{split}
 \eeqlb
Combining $(\ref{t5.1.5})$ and $(\ref{t5.1.7b})$, we obtain
 \beqnn
 \E\big(T(r)^n\big)&\leq&cnr^Q\int_{I^{n-1}}\P\left\{\|X(t^1)\|<r,
 \cdots,\|X(t^{n-1})\|<r\right\}dt^1\cdots dt^{n-1}\cr
 &=&cnr^Q \E\big(T(r)^{n-1}\big).
 \eeqnn
Hence the inequality $(\ref{l5.1.1})$ follows from this
and induction. \qed

Let $0<b<1/c_{_{3, 16}}$. Then by $(\ref{l5.1.1})$ we have
 \beqn
 \label{l5.1.2}
 \E\left(\exp\left(br^{-Q}T(r)\right)\right)\leq
 \sum_{n=0}^\infty(c_{_{3, 16}}b)^n<\infty.
 \eeqn
This and the exponential Chebychev's inequality imply that
for any constant $0<b<1/c_{_{3, 16}}$,
\beqn
\label{l5.1.2b}
\P\big\{T(r) \ge r^Q u\big\} \le \frac {e^{-bu}} {1- c_{_{3, 16}}b}
\eeqn
for all $u > 0$.

The following is a law of the iterated logarithm for the sojourn
measure of $X$.
\begin{proposition} \label{p5.2}
For every $\tau\in I$, we have
\beqlb
\label{4.8*}
\limsup_{r\rightarrow0}\frac{T_{X(\tau)}(r)}{\varphi_1(r)}
\leq c_{_{3, 16}}, \quad a.s.
\eeqlb
\end{proposition}

\noindent\textbf{Proof.} Since $\{X(t), t\in\N\}$ has stationary
increments, it is sufficient to consider $\tau = 0$. Then
$(\ref{4.8*})$ follows from (\ref{l5.1.2b}), the
Borel-Cantelli lemma and a monotonicity argument
in a standard way. \qed

\noindent\textbf{Proof of Theorem $\ref{t5.1}$.} We can prove this
theorem by using Lemma $\ref{l3.1}$ and Proposition $\ref{p5.2}$,
in the same way as that of Theorem 4.1 in Xiao (1996).
 \qed

\noindent\textbf{Proof of Theorem \ref{t1.1}.} It follows
immediately from Theorems $\ref{t4.1}$ and $\ref{t5.1}$.
\qed

\bibliographystyle{plain}

\begin{thebibliography}{1234}

\bibitem{Adler81}
Adler, R. J. (1981), {\it The Geometry of Random Fields}. Wiley, New
York.



\bibitem{Anderson55}
Anderson, T. W. (1955), The integral of a symmetric unimodal
function over a symmetric convex set and some probability
inequalities. {\it Proc. Amer. Math. Soc.} {\bf 6}, 170-176.

\bibitem{AyacheXiao05}
Ayache, A. and Xiao, Y. (2005), Asymptotic growth properties and
Hausdorff dimension of fractional Brownian sheets. {\it J. Fourier
Anal. Appl.} {\bf 11}, 407--439.

\bibitem{B-M1}
Baraka, D. and Mountford, T. (2008), A law of iterated logarithm for
fractional Brownian motions. {\it S\'eminaire de Probabilit\'es XLI.
Lecture Notes in Math.} {\bf 1934}, 161--179, Springer, Berlin.

\bibitem{B&M2}
Baraka, D. and Mountford, T. (2011), The exact Hausdorff measure of
the zero set of fractional Brownian motion. {\it J. Theor. Probab.}
{\bf 24}, 271-�293.


\bibitem{B&J&R1997}
Benassi, A., Jaffard, S. and Roux, D. (1997), Elliptic Gaussian
random Processes. \emph{Rev. Mat. Iberoamericana} \textbf{13},
19--90.

\bibitem{BergForst75}
Berg, C. and Forst, G. (1975), {\it Potential Theory on Locally
Compact Abelian Groups}. Sringer-Verlag, New York-Heidelberg.


\bibitem{Berman73}
Berman, S. M. (1973), Local nondeterminism and local times of
Gaussian processes.  {\it Indiana Univ. Math. J.} {\bf 23}, 69--94.

\bibitem{Berman88}
Berman, S. M. (1988),  Spectral conditions for local nondeterminism.
{\it Stochastic Process. Appl.} {\bf 27}, 73--84.

\bibitem{BLX09}
Bierm\'e, H., Lacaux,  C. and Xiao, Y. (2009), Hitting probabilities and
the Hausdorff dimension of the inverse images of anisotropic
Gaussian random fields. {\it Bull. London Math. Soc.} {\bf 41},
253--273.




\bibitem{Fal90}
Falconer, K. J. (1990), {\it Fractal Geometry -- Mathematical
Foundations and Applications}. Wiley \& Sons, New York.

\bibitem{Istas05}
Istas, J. (2005), Spherical and hyperbolic fractional Brownian motion.
{\it Elec. Comm. Probab.} {\bf 10}, 254--262.

\bibitem{Kahane85}
Kahane, J.-P. (1985), {\it Some Random Series of Functions}. 2nd
edition, Cambridge University Press, Cambridge.

\bibitem{L77}
Lo\'{e}ve, L. (1977), {\it Probability Theory I.} Springer-Verlag,
New York.

\bibitem{L&X2010}
Luan, N. and Xiao, Y. (2010), Chung's law of the iterated logarithm
for anisotropic Gaussian random fields. {\it Statist. Probab. Lett.}
{\bf 80}, 1886--1895.



\bibitem{MWX10}
Meerschaert, M. M., Wang, W. and Xiao,  Y. (2011), Fernique-type
inequalities and moduli of continuity of anisotropic Gaussian random
fields. {\it Trans. Amer. Math. Soc.}, to appear.


\bibitem{M&R87}
Monrad, D. and Rootz\'{e}n, H. (1995), Small values of Gaussian
processes and functional laws of the iterated logarithm. {\it Proba.
Theory Relat. Fields} \textbf{101}, 173--192.

\bibitem{NuaViens08}
Nualart, E. and Viens, F. (2009), The fractional stochastic heat
equation on the circle: time regularity and potential theory. {\it
Stoch. Process. Appl.} {\bf 119},  1505--1540.

\bibitem{Pitt75}
Pitt, L. D. (1975), Stationary Gaussian Markov fields on $\R^d$ with
a deterministic component. {\it J. Multivar. Anal.} {\bf 5},
300--311.

\bibitem{Pitt78}
Pitt, L. D. (1978), Local times for Gaussian vector fields. {\it
Indiana Univ. Math. J.} {\bf 27}, 309--330.

\bibitem {RT}
Rogers, C. A. and Taylor, S. J. (1961), Functions continuous and
singular with respect to a Hausdorff measure. {\it Mathematika} {\bf
8}, 1--31.

\bibitem{[T93]}
Talagrand, M. (1993), New Gaussian estimates for enlarged balls.
\textit{Geom. Funct. Anal.} \textbf{3}, 502--526.

\bibitem{Talagrangd95}
Talagrand, M. (1995), Hausdorff measure of  trajectories of
multiparameter fractional Brownian motion.  {\it Ann. Probab.} {\bf
23}, 767--775.


\bibitem{Talagrand98}
Talagrand, M. (1998), Multiple points of trajectories of
multiparameter fractional Brownian motion. {\it Probab. Theory
Relat. Fields} {\bf 112}, 545--563.


\bibitem{TTV04}
Tindel, S., Tudor, C. A. and Viens, F. (2004), Sharp Gaussian
regularity on the circle, and applications to the fractional stochastic
heat equation. {\it J. Funct. Anal.} {\bf 217}, 280--313.



\bibitem{WX09}
Wu, D.  and Xiao, Y. (2009), Uniform Hausdorff dimension results for
Gaussian random fields. {\it Sci. in China, Ser. A} {\bf 52},
1478--1496.

\bibitem{WX10}
Wu, D.  and Xiao, Y. (2011), On local times of anisotropic Gaussian
random fields. {\it Comm. Stoch. Anal.} {\bf 5}, 15--39.

\bibitem{Xiao96}
Xiao, Y. (1996), Hausdorff measure of the sample paths of Gaussian
random fields. {\it Osaka J. Math.} {\bf 33}, 895--913.

\bibitem{Xiao97a}
Xiao, Y. (1997a), Hausdorff dimension of the graph of fractional
Brownian motion. {\it Math. Proc. Camb. Philos. Soc.} {\bf 122},
565--576.

\bibitem{Xiao97b}
Xiao, Y. (1997b), H\"older conditions for the local times and the
Hausdorff measure of the level sets of Gaussian random fields. {\it
Probab. Theory Relat. Fields} {\bf 109}, 129--157.



\bibitem{Xiao07}
Xiao, Y. (2007),  Strong local nondeterminism of Gaussian random
fields and its applications. In: {\it Asymptotic Theory in
Probability and Statistics with Applications}, (T.-L. Lai, Q.-M.
Shao and L. Qian, editors), pp. 136--176, Higher Education Press,
Beijing.


\bibitem{Xiao09}
Xiao, Y. (2009),  Sample path properties of anisotropic Gaussian
random fields. In: {\it A Minicourse on Stochastic Partial
Differential Equations,} (D. Khoshnevisan and F. Rassoul-Agha,
editors), Lecture Notes in Math. {\bf 1962}, pp. 145--212, Springer,
New York.

\bibitem{XueXiao09}
Xue, Y. and Xiao, Y. (2011), Fractal and smoothness properties of
anisotropic Gaussian models. {\it Frontiers Math. China}, to appear.

\bibitem{Yaglom57}
Yaglom, A. M. (1957),  Some classes of random fields in
$n$-dimensional space, related to stationary random processes. {\it
Th. Probab. Appl.} {\bf 2}, 273--320.
\end{thebibliography}
\begin{small}

\end{small}

\bigskip

\noindent \textsc{Nana Luan}.\  School of Insurance and Economics,
University
of International Business and Economics, Beijing 100029, China\\
        E-mail: \texttt{luannana318@gmail.com}
        \\

\noindent \textsc{Yimin Xiao}.\
        Department of Statistics and Probability, A-413 Wells
        Hall, Michigan State University,
        East Lansing, MI 48824, U.S.A.\\
        E-mail: \texttt{ xiao@stt.msu.edu}\\
        URL: \texttt{http://www.stt.msu.edu/\~{}xiaoyimi}

\end{document}